\documentclass[12pt, a4paper]{article}
\usepackage[pagewise]{lineno}
\usepackage{graphicx}
\usepackage{amssymb}
\usepackage{latexsym, bm}
\usepackage{multicol}
\usepackage{indentfirst}
\usepackage{amssymb,amsfonts}
\usepackage{amsmath}
\usepackage{cases}
\usepackage[noadjust]{cite}
\usepackage[colorlinks,citecolor=red]{hyperref}
\newtheorem{theorem}{Theorem}
\newtheorem{lemma}{Lemma}

\newtheorem{definition}{Definition}

\newtheorem{remark}{Remark}

\usepackage{amssymb}
\numberwithin{equation}{section}
\numberwithin{theorem}{section}
\numberwithin{remark}{section}
\numberwithin{definition}{section}
\numberwithin{lemma}{section}
\title{Regularity criteria via one directional derivative of the velocity in anisotropic Lebesgue spaces to the 3D Navier-Stokes equations}
\author{Maria Alessandra Ragusa\footnote{maragusa@dmi.unict.it}\\
{\small Department of Mathematics, University of Catania,}\\
{\small Viale Andrea Doria No. 6, Catania 95128, Italy}\\
{\small RUDN University,}\\
{\small6 Miklukho -Maklay St, Moscow, 117198, Russia}
\and
Fan Wu\footnote{wufan0319@yeah.net}\\
{\small School of Mathematics and Statistics, Hunan Normal University, }\\
{\small Changsha, Hunan 410081, China }}

\date{}

\begin{document}
\maketitle
{\bf Abstract:}
In this paper, we consider the regularity criterion for 3D incompressible Navier-Stokes equations in terms of one directional derivative of the velocity in anisotropic Lebesgue spaces.  More precisely, it is proved that u becomes a regular solution if the $\partial_3u$ satisfies
 $$\int^{T}_{0}\frac{\left\|\left\|\left\|\partial_{3}u(t)\right\|_{L^{p}_{x_1}}\right\|_{L^{q}_{x_2}}\right\|^{\beta}_{L^{r}_{x_3}}}
 {1+\ln\left(\|\partial_3u\left(t\right)\|_{L^2}+e\right)}dt<\infty,$$
$\text { where } \frac{2}{\beta}+\frac{1}{p}+\frac{1}{q}+\frac{1}{r}=1 \text { and } 2 < p, q, r \leq \infty, 1-\left(\frac{1}{p}+\frac{1}{q}+\frac{1}{r}\right) \geq 0
$.

\medskip

\section{Introduction}
In this paper, we are concerned with regularity criteria for the weak solutions to the incompressible Navier-Stokes
equations (NSE) in $\mathbb{R}^3$ :
\begin{equation}
\begin{gathered}\label{1.1}
\left\{
             \begin{array}{lr}
             \partial_t u+(u\cdot\nabla)u-\Delta u+\nabla p=0,& \\
             \nabla\cdot u=0,&\\
             u(x,0)=u_{0}(x),&
\end{array}
\right.
\end{gathered}
\end{equation}
where $u=\left(u_{1}\left(t, x\right), u_{2}\left(t, x\right), u_{3}\left(t, x\right)\right)$ and $p$ denote the unknown velocity field and pressure of the fluid respectively, and $u_{0}$ is the prescribed initial data satisfying the compatibility condition $\nabla\cdot u_{0}$ = 0, and
$$
\partial_{t} u=\frac{\partial u}{\partial t}, \quad \partial_{i}=\frac{\partial}{\partial x_{i}}, \quad(u \cdot \nabla)=\sum_{i=1}^{3} u_{i} \partial_{i}.
$$

It is well known that the weak solution of Navier-Stokes equations \eqref{1.1} is unique and regular in two dimensions. However, in three dimensions, the regularity problem of weak solutions of Navier-Stokes equations is an outstanding open problem in mathematical mechanics. The classical Prodi-Serrin conditions \cite{GP,SJ,ELSGA} say that if
\begin{equation}\label{1.2}
u \in L^{q}\left(0, T ; L^{p}\left(\mathbb{R}^{3}\right)\right) \quad \text{with} \quad\frac{2}{q}+\frac{3}{p}=1 \quad \text{and}\quad 3\leq p\leq\infty,
\end{equation}
then the weak solution $u$ is regular on $(0,T]$. The Prodi-Serrin conditions \eqref{1.2} were later generalized by Beir$\tilde{a}$oda Veiga \cite{HBV} to be
\begin{equation}\label{1.3}
\nabla u \in L^{q}\left(0, T ; L^{p}\left(\mathbb{R}^{3}\right)\right)\quad \text{with} \quad\frac{2}{q}+\frac{3}{p}=2 \quad \text{and}\quad \frac{3}{2}\leq p\leq\infty.
\end{equation}
In the past decades, many refinements of \eqref{1.2} and \eqref{1.3} appeared. Penel and Pokorny \cite{PM} first established a regularity criterion for the 3D Navier-Stokes equations
only in terms of one directional derivative of the velocity, more precisely, they proved
that if
\begin{equation}\label{1.4}
\partial_{3} u \in L^{q}\left(0, T ; L^{p}\left(\mathbb{R}^{3}\right)\right) \quad \text{with} \quad\frac{2}{q}+\frac{3}{p}=\frac{3}{2} \quad \text{and}\quad 2 \leq p \leq \infty,
\end{equation}
then the solution is smooth. Kukavica and Ziane \cite{KIZM} extended the above condition to the condition
\begin{equation}\label{1.5}
\partial_{3} u \in L^{q}\left(0, T ; L^{p}\left(\mathbb{R}^{3}\right)\right) \quad \text{with} \quad\frac{2}{q}+\frac{3}{p}=2 \quad \text{and}\quad \frac{9}{4} \leq p \leq 3.
\end{equation}
Recently, Cao \cite{CCS} and Zhang et al. \cite{ZZJ,ZYZ} have made various improvements to the indicators in condition \eqref{1.5}. In \cite{LQ},  Liu proved the regularity criterion via one directional derivative of the velocity (i.e. $\partial_3u$) in Morrey-Campanato spaces. For
readers interested in this topic for partial components and related arguments, please refer to \cite{FDY,SK,SGRM, GR2012, GR2016, PR, ZZF,GG} for recent progresses.

Motivated by papers cited above, we shall investigate regularity criterion for the weak solutions to the Cauchy problem of the 3D Navier-Stokes equations \eqref{1.1} in term of the one directional derivative of the velocity on framework of the anisotropic Lebesgue spaces in this
work. Before stating our main theorem, we shall first recall the definitions of the anisotropic Lebesgue spaces and the weak solutions to \eqref{1.1}.
\begin{definition}
Let $1 \leq p, q, r \leq \infty$, we say that a function $f$ belongs to anisotropic Lebesgue spaces $L^p\left(\mathbb{R}_{x_{1}} ; L^{q}\left(\mathbb{R}_{x_{2}} ; L^{r}\left(\mathbb{R}_{x_{3}}\right)\right)\right)$ if $f$ is measurable on $\mathbb{R}^3$
and the following norm is finite:
$$
\left\|\left\| \left\| f\right\|_{L_{x_{1}}^{p}}\right\|_{L_{x_{2}}^{q}} \right\|_{L_{x_{3}}^{r}} :=\left(\int_{\mathbb{R}}\left(\int_{\mathbb{R}}\left(\int_{\mathbb{R}}\left|f\left(x_{1}, x_{2}, x_{3}\right)\right|^{p} \mathrm{d} x_{1}\right)^{\frac{q}{p}} \mathrm{d} x_{2}\right)^{\frac{r}{q}} \mathrm{d} x_{3}\right)^{\frac{1}{r}}
$$
with the usual change as $p=\infty$ or $q=\infty$ or $r=\infty$.
\end{definition}
\begin{definition}
Let $u_0\in L^2(\mathbb{R}^3)$ with $\nabla\cdot u_0 = 0$, $T > 0$. A measurable $\mathbb{R}^3$-valued function $u$ defined in $[0, T ] \times \mathbb{R}^3$ is said to be a weak solution to \eqref{1.1} if\\
1. $u \in L^{\infty}\left(0, T ; L^{2}\left(\mathbb{R}^{3}\right)\right) \cap L^{2}\left(0, T ; H^{1}\left(\mathbb{R}^{3}\right)\right)$; \\
2. $\eqref{1.1}_1$ holds in the sense of distributions, that is,
$$
\int_{0}^{t}\left(u, \partial_{t} \phi+(u \cdot \nabla)\phi\right) d s+\left(u_{0}, \phi(0)\right)=\int_{0}^{t}\left(\nabla u\left(s\right), \nabla \phi\left(s\right)\right) d s,
$$
for all $\phi\in C^{\infty}_{c}\left([0, T ) \times \mathbb{R}^3\right) $ with $\nabla\cdot\phi=0$, where $(\cdot, \cdot)$ is the scalar product in $L^2(\mathbb{R}^3)$;\\
3.  the energy inequality, that is,
$$
\|u(t)\|_{L^{2}}^{2}+2 \int_{0}^{t}\|\nabla u(s)\|_{L^{2}}^{2} d s \leq\left\|u_{0}\right\|_{L^{2}}^{2}, \quad 0 \leq t \leq T.
$$
\end{definition}
Now, our main result reads:
\begin{theorem}\label{th1}
 Suppose that $u_{0}\in H^{1}(\mathbb{R}^3)$ and $\nabla\cdot u_{0}=0$ in distribution. Let $u$ be the Leray-Hopf weak solution of \eqref{1.1} on $(0,T]$. If $\partial_3u$  satisfies
\begin{equation}\label{1.6}
\int^{T}_{0}\frac{\left\|\left\|\left\|\partial_{3}u(t)\right\|_{L^{p}_{x_1}}\right\|_{L^{q}_{x_2}}\right\|^{\beta}_{L^{r}_{x_3}}}
 {1+\ln\left(\|\partial_3u\left(t\right)\|_{L^2}+e\right)}dt<\infty,
\end{equation}
where $\frac{2}{\beta}+\frac{1}{p}+\frac{1}{q}+\frac{1}{r}=1 $ and $ 2 < p, q, r \leq \infty, 1-\left(\frac{1}{p}+\frac{1}{q}+\frac{1}{r}\right) \geq 0$,
then the weak solution $u$ is smooth on interval $(0,T]$.
\end{theorem}
\begin{remark}
The result of Theorem \ref{th1} make a further step forward in the understanding of regularity criteria for the Navier-Stokes equations. Due to
$$
\frac{\left\|\left\|\left\|\partial_{3}u(t)\right\|_{L^{p}_{x_1}}\right\|_{L^{q}_{x_2}}\right\|^{\beta}_{L^{r}_{x_3}}}
 {1+\ln\left(\|\partial_3u\left(t\right)\|_{L^2}+e\right)}
\leq\left\|\left\|\left\|\partial_{3}u(t)\right\|_{L^{p}_{x_1}}\right\|_{L^{q}_{x_2}}\right\|^{\beta}_{L^{r}_{x_3}},
$$
we can easily get regularity condition
\begin{equation}\label{1.7}
\int^{T}_{0}\left\|\left\|\left\|\partial_{3}u(t)\right\|_{L^{p}_{x_1}}\right\|_{L^{q}_{x_2}}\right\|^{\beta}_{L^{r}_{x_3}}
 dt<\infty,
 \end{equation}
where $\frac{2}{\beta}+\frac{1}{p}+\frac{1}{q}+\frac{1}{r}=1 $ and $ 2 < p, q, r \leq \infty, 1-\left(\frac{1}{p}+\frac{1}{q}+\frac{1}{r}\right) \geq 0$.
\end{remark}
\begin{remark}
It seems a difficult problem to prove the integrable index  $\frac{2}{\beta}+\frac{1}{p}+\frac{1}{q}+\frac{1}{r}=2 $ in \eqref{1.6}. We
hope we can overcome this problem in the near future.
\end{remark}

Before ending this section, we state the following two lemmas, which will be used in
the proof of our main results.
\begin{lemma}\cite{LIU}\label{le1}
There exists a positive constant C such that
\begin{equation}\label{1.8}
\left\| \left\| \left\|f\right\|_{L_{x_{1}}^{\frac{2 p}{p-2}}}\right\|_{L_{x_{2}}^{\frac{2 q}{q-2}}} \right\|_{L_{x_{3}}^{\frac{2 r}{r-2}}}\leq C\left\|\partial_{1} f\right\|_{L^{2}}^{\frac{1}{p}}\left\|\partial_{2} f\right\|_{L^{2}}^{\frac{1}{q}}\left\|\partial_{3} f\right\|_{L^{2}}^{\frac{1}{r}}\|f\|_{L^{2}}^{1-\left(\frac{1}{p}+\frac{1}{q}+\frac{1}{r}\right)},
\end{equation}
for every $f \in C_{0}^{\infty}\left(\mathbb{R}^{3}\right)$ where $2 < p, q, r \leq \infty, 1-\left(\frac{1}{p}+\frac{1}{q}+\frac{1}{r}\right) \geq 0$.
\end{lemma}
\begin{lemma}\label{le2}
\cite{CW} Let $\mu,\theta,\lambda$ and $\kappa$ be four numbers satisfying
\begin{displaymath}
1\leq\mu,\theta,\lambda,\kappa<\infty,\quad\frac{1}{\theta}+\frac{1}{\lambda}+\frac{1}{\kappa}>1
\quad and\quad 1+\frac{3}{\mu}=\frac{1}{\theta}+\frac{1}{\lambda}+\frac{1}{\kappa}.
\end{displaymath}
Assume that $\varphi(x)=\varphi(x_{1},x_{2},x_{3})$ with $\partial_{1}\varphi\in L^{\theta}(\mathbb{R}^{3}),\partial_{2}\varphi\in L^{\lambda}(\mathbb{R}^{3})$ and $\partial_{3}\varphi\in L^{\kappa}(\mathbb{R}^{3})$.
Then, there exists a constant $C=C(\theta,\lambda,\kappa)$ such that
\begin{equation}\label{1.9}
\|\varphi\|_{L^{\mu}}\leq C\|\partial_{1}\varphi\|^{\frac{1}{3}}_{L^{\theta}}
\|\partial_{2}\varphi\|^{\frac{1}{3}}_{L^{\lambda}}
\|\partial_{3}\varphi\|^{\frac{1}{3}}_{L^{\kappa}}.
\end{equation}
Especially, when $\theta=\lambda=2$ and $1\leq\kappa<\infty$, there exists a constant $C=C(\kappa)$ such that
\begin{equation}\label{1.10}
\|\varphi\|_{L^{3\kappa}}\leq C\|\partial_{1}\varphi\|^{\frac{1}{3}}_{L^{2}}
\|\partial_{2}\varphi\|^{\frac{1}{3}}_{L^{2}}
\|\partial_{3}\varphi\|^{\frac{1}{3}}_{L^{\kappa}},
\end{equation}
which holds for any $\varphi$ with $\partial_{1}\varphi\in L^{2}(\mathbb{R}^{3}), \partial_{2}\varphi\in L^{2}(\mathbb{R}^{3})$ and $\partial_{3}\varphi\in L^{\kappa}(\mathbb{R}^{3})$.
\end{lemma}
Throughout this paper, we denote by $C$ a  positive constant,
which may depend on the initial data $u_0$ and $T$ , and its value may change from line to
line. The norms of the usual Lebesgue spaces $L^p(\mathbb{R}^3)$ are denoted by $L^p$ , while
the directional derivatives of a function $f$ are denoted by $\partial_{i} f=\frac{\partial f}{\partial x_{i}}(i=1,2,3)$.

We present the proof of Theorems \ref{th1} in the next section.
\section{The proof of Theorem 1.1}
This section is devoted to the proof of Theorem \ref{th1}. The proof is based on the establishment of a priori estimates under condition \eqref{1.6}.

Firstly, we give the basic energy estimate of equations \eqref{1.1}. Multiplying the first equation of \eqref{1.1} with $u$ and integrating them in ${\mathbb{R}}^3$, we have
\begin{equation}\label{2.1}
\frac{1}{2}\frac{d}{dt}\|u\|^{2}_{L^2}+\|\nabla u\|^{2}_{L^2}=0.
\end{equation}
Integrating from 0 to $t$ for the above equality, we have
\begin{equation}\label{2.2}
u\in L^{\infty}\big(0,T;L^{2}(\mathbb{R}^{3})\big)\cap L^{2}\big(0,T;H^{1}(\mathbb{R}^{3})\big).
\end{equation}

Differentiating the first equation in $(1.1)$ about space variable $x_{3}$, then multiplying the resulting equation by $\partial_{3}u$, and integrating it to get
\begin{equation}\label{2.3}
\begin{split}
\frac{1}{2}\frac{d}{dt}\|\partial_{3} u\|^{2}_{L^{2}}+\|\nabla \partial_{3}u\|^{2}_{L^{2}}&=-\int_{\mathbb{R}^3}\partial_{3}(u\cdot\nabla u)\cdot\partial_{3} u dx=-\int_{\mathbb{R}^3}\partial_{3}u\cdot\nabla u\cdot\partial_{3} u dx.\\
\end{split}
\end{equation}
By H\"older's inequality, Young's inequality and Lemma \ref{le1}, we can estimate the right-hand side of equality \eqref{2.3} as follows:
\begin{equation}
\begin{split}\label{2.4}
&-\int_{\mathbb{R}^3}\partial_{3}u\cdot\nabla u\cdot\partial_{3} u dx\leq\left|\int_{\mathbb{R}^3}\partial_{3}u\cdot\nabla u\cdot\partial_{3} u dx\right|\\
&\leq\left\|\left\|\left\|\partial_{3}u\right\|_{L^{p}_{x_1}}\right\|_{L^{q}_{x_2}}\right\|_{L^{r}_{x_3}}
\left\| \left\| \left\|\partial_3u\right\|_{L_{x_{1}}^{\frac{2 p}{p-2}}}\right\|_{L_{x_{2}}^{\frac{2 q}{q-2}}} \right\|_{L_{x_{3}}^{\frac{2 r}{r-2}}}\left\|\nabla u\right\|_{L^2}\\
&\leq C\left\|\left\|\left\|\partial_{3}u\right\|_{L^{p}_{x_1}}\right\|_{L^{q}_{x_2}}\right\|_{L^{r}_{x_3}}
\left\|\partial_{1}\partial_{3}u \right\|_{L^{2}}^{\frac{1}{p}}\left\|\partial_{2} \partial_{3}u\right\|_{L^{2}}^{\frac{1}{q}}\left\|\partial_{3} \partial_{3}u\right\|_{L^{2}}^{\frac{1}{r}}\|\partial_{3}u\|_{L^{2}}^{1-\left(\frac{1}{p}+\frac{1}{q}+\frac{1}{r}\right)}\|\nabla u\|_{L^{2}}\\
&\leq C\left\|\left\|\left\|\partial_{3}u\right\|_{L^{p}_{x_1}}\right\|_{L^{q}_{x_2}}\right\|_{L^{r}_{x_3}}
\left\|\nabla\partial_{3}u \right\|_{L^{2}}^{\frac{1}{p}+\frac{1}{q}+\frac{1}{r}}\|\partial_{3}u\|_{L^{2}}^{1-\left(\frac{1}{p}+\frac{1}{q}+\frac{1}{r}\right)}\|\nabla u\|_{L^{2}}\\
&\leq C\left(\|\nabla\partial_{3}u\|^{2}_{L^{2}}\right)^{\frac{\frac{1}{p}+\frac{1}{q}+\frac{1}{r}}{2}}
\left(\left\|\left\|\left\|\partial_{3}u\right\|_{L^{p}_{x_1}}\right\|_{L^{q}_{x_2}}\right\|^{\frac{2}{2-\left(\frac{1}{p}+\frac{1}{q}+\frac{1}{r}\right)}}_{L^{r}_{x_3}}
\|\partial_{3}u\|^{2\cdot\frac{1-\left(\frac{1}{p}+\frac{1}{q}+\frac{1}{r}\right)}{2-\left(\frac{1}{p}+\frac{1}{q}+\frac{1}{r}\right)}}_{L^{2}}\|\nabla u\|^{\frac{2}{2-\left(\frac{1}{p}+\frac{1}{q}+\frac{1}{r}\right)}}_{L^{2}}\right)^{\frac{2-\left(\frac{1}{p}+\frac{1}{q}+\frac{1}{r}\right)}{2}}\\
&\leq \epsilon\|\nabla\partial_{3}u\|^{2}_{L^{2}}+C\left(\left\|\left\|\left\|\partial_{3}u\right\|_{L^{p}_{x_1}}\right\|_{L^{q}_{x_2}}\right\|^{\frac{2}{2-\left(\frac{1}{p}+\frac{1}{q}+\frac{1}{r}\right)}}_{L^{r}_{x_3}}
\|\partial_{3}u\|^{2\cdot\frac{1-\left(\frac{1}{p}+\frac{1}{q}+\frac{1}{r}\right)}{2-\left(\frac{1}{p}+\frac{1}{q}+\frac{1}{r}\right)}}_{L^{2}}\|\nabla u\|^{\frac{2}{2-\left(\frac{1}{p}+\frac{1}{q}+\frac{1}{r}\right)}}_{L^{2}}\right)\\
&\leq \epsilon\|\nabla\partial_{3}u\|^{2}_{L^{2}}+C\left(\left\|\left\|\left\|\partial_{3}u\right\|_{L^{p}_{x_1}}\right\|_{L^{q}_{x_2}}\right\|_{L^{r}_{x_3}}^{\frac{2}{1-\left(\frac{1}{p}+\frac{1}{q}+\frac{1}{r}\right)}}\right)^{\frac{1-\left(\frac{1}{p}+\frac{1}{q}+\frac{1}{r}\right)}{2-\left(\frac{1}{p}+\frac{1}{q}+\frac{1}{r}\right)}}
\left(\|\nabla u\|^{2}_{L^{2}}\right)^{\frac{1}{2-\left(\frac{1}{p}+\frac{1}{q}+\frac{1}{r}\right)}}
\|\partial_{3}u\|^{2\cdot\frac{1-\left(\frac{1}{p}+\frac{1}{q}+\frac{1}{r}\right)}{2-\left(\frac{1}{p}+\frac{1}{q}+\frac{1}{r}\right)}}_{L^{2}}\\
&\leq \epsilon\|\nabla\partial_{3}u\|^{2}_{L^{2}}+C\left(1+\|\partial_{3}u\|^{2}_{L^{2}}\right)\left(
\left\|\left\|\left\|\partial_{3}u\right\|_{L^{p}_{x_1}}\right\|_{L^{q}_{x_2}}\right\|_{L^{r}_{x_3}}^{\frac{2}{1-\left(\frac{1}{p}+\frac{1}{q}+\frac{1}{r}\right)}}+\|\nabla u\|^{2}_{L^{2}}\right).
\end{split}
\end{equation}
Combining \eqref{2.3} and \eqref{2.4} together, then we get
\begin{equation}
\begin{split}\label{2.5}
&\frac{d}{dt}(e+\|\partial_{3} u\|^{2}_{L^{2}})+\|\nabla \partial_{3}u\|^{2}_{L^{2}}\\
\leq &C\left(e+\|\partial_{3}u\|^{2}_{L^{2}}\right)\left(
\left\|\left\|\left\|\partial_{3}u\right\|_{L^{p}_{x_1}}\right\|_{L^{q}_{x_2}}\right\|_{L^{r}_{x_3}}^{\frac{2}{1-\left(\frac{1}{p}+\frac{1}{q}+\frac{1}{r}\right)}}+\|\nabla u\|^{2}_{L^{2}}\right)\\
\leq &C\left(e+\|\partial_{3}u\|^{2}_{L^{2}}\right)\left(\frac{\left\|\left\|\left\|\partial_{3}u(t)\right\|_{L^{p}_{x_1}}\right\|_{L^{q}_{x_2}}\right\|^{\beta}_{L^{r}_{x_3}}}
 {e+\ln\left(\|\partial_3u\left(t\right)\|_{L^2}+e\right)}+\|\nabla u\|^{2}_{L^{2}}\right)\left(1+\ln\left(\|\partial_3u\left(t\right)\|^2_{L^2}+e\right)\right).\\
\end{split}
\end{equation}
Let
$$
F(t)=\ln\left(\|\partial_3u\left(t\right)\|^2_{L^2}+e\right),
$$
we obtain
\begin{equation}
\begin{split}\label{2.6}
\frac{d}{dt}\left(1+F\left(t\right)\right)+\|\nabla \partial_{3}u\|^{2}_{L^{2}}
\leq C\left(\frac{\left\|\left\|\left\|\partial_{3}u(t)\right\|_{L^{p}_{x_1}}\right\|_{L^{q}_{x_2}}\right\|^{\beta}_{L^{r}_{x_3}}}
 {1+\ln\left(\|\partial_3u\left(t\right)\|_{L^2}+e\right)}+\|\nabla u\|^{2}_{L^{2}}\right)\left(1+F(t)\right).
\end{split}
\end{equation}
Applying the Gronwall inequality to \eqref{2.6} yields that
 \begin{equation}
\begin{split}\label{3.16}
\sup\limits_{0\leq t\leq T}\ln F(t)&\leq \left(1+\ln F(0)\right)\exp \int^{T}_{0}C\left(\frac{\left\|\left\|\left\|\partial_{3}u(t)\right\|_{L^{p}_{x_1}}\right\|_{L^{q}_{x_2}}\right\|^{\beta}_{L^{r}_{x_3}}}
 {1+\ln\left(\|\partial_3u\left(t\right)\|_{L^2}+e\right)}+\|\nabla u\|^{2}_{L^{2}}\right)dt,
\end{split}
\end{equation}
which implies that
\begin{equation}
\begin{split}\label{2.8}
\sup\limits_{0\leq t\leq T}\|\partial_{3} u(t)\|^{2}_{L^{2}}+\int^{T}_{0}\|\nabla\partial_{3} u(t)\|^{2}_{L^{2}}dt&\leq C.
\end{split}
\end{equation}
Next, multiplying the first equation of \eqref{1.1} by $\Delta u$, integrating over $\mathbb{R}^3$, we get
 \begin{equation}
 \begin{split}\label{2.9}
\frac{1}{2}\frac{d}{dt}\|\nabla u\|^{2}_{L^{2}}+\|\Delta u\|^{2}_{L^{2}}&=\int_{\mathbb{R}^3}u\cdot\nabla u\cdot\Delta u dx.
\end{split}
\end{equation}
For the term of right-hand side of \eqref{2.9}, by using the interpolation inequality and Lemma \ref{le2} (with $\kappa=2$ in \eqref{1.10}), it is not difficult to see that
\begin{equation}
\begin{split}\label{2.10}
\int_{\mathbb{R}^3}u\cdot\nabla u\cdot\Delta u dx&\leq C\|\nabla u\|^{3}_{L^{3}}\leq C\|\nabla u\|^{\frac{3}{2}}_{L^{2}}\|\nabla u\|^{\frac{3}{2}}_{L^{6}}\\
&\leq C\|\nabla u\|^{\frac{3}{2}}_{L^{2}}\|\nabla\partial_{1}u\|^{\frac{1}{2}}_{L^{2}}
\|\nabla\partial_{2}u\|^{\frac{1}{2}}_{L^{2}}\|\nabla\partial_{3}u\|^{\frac{1}{2}}_{L^{2}}\\
&\leq C\|\nabla u\|^{\frac{3}{2}}_{L^{2}}\|\nabla^{2}u\|_{L^{2}}
\|\nabla\partial_{3}u\|^{\frac{1}{2}}_{L^{2}}\\
&=C(\|\nabla^{2}u\|^{2}_{L^{2}})^{\frac{1}{2}}(\|\nabla u\|^{3}_{L^{2}}\|\nabla\partial_{3}u\|_{L^{2}})^{\frac{1}{2}}\\
&\leq\frac{1}{4}\|\Delta u\|^{2}_{L^{2}}+C\|\nabla u\|^{3}_{L^{2}}\|\nabla\partial_{3}u\|_{L^{2}}\\
&\leq\frac{1}{4}\|\Delta u\|^{2}_{L^{2}}+C\|\nabla u\|^{2}_{L^{2}}(\|\nabla u\|^2_{L^{2}}+\|\nabla\partial_{3}u\|^2_{L^{2}}).
\end{split}
\end{equation}
Substituting \eqref{2.10} into \eqref{2.9}, we obtain
 \begin{equation}
 \begin{split}\label{2.11}
\frac{d}{dt}\|\nabla u\|^{2}_{L^{2}}+\|\Delta u\|^{2}_{L^{2}}
&\leq C\|\nabla u\|^{2}_{L^{2}}(\|\nabla\partial_{3}u\|^{2}_{L^{2}}+\|\nabla u\|^{2}_{L^{2}}).\\
\end{split}
\end{equation}
Applying the Gronwall inequality yields that
\begin{equation}
\begin{split}
\sup\limits_{0\leq t\leq T}&\|\nabla u\|^{2}_{L^{2}}\leq \|\nabla u_{0}\|^{2}_{L^{2}}\exp C\left\{\int^{T}_{0}\left(\|\nabla\partial_{3}u\|^{2}_{L^{2}}+\|\nabla u\|^{2}_{L^{2}}\right)d\tau\right\},
\end{split}
\end{equation}
which gives that
$$u\in L^{\infty}\left(0,T;H^{1}\left(\mathbb{R}^3\right)\right)\cap L^{2}\left(0,T;H^{2}\left(\mathbb{R}^3\right)\right).$$
This completes the proof of Theorem \ref{th1}.
\section*{Acknowledgments}
\par The first author is partially supported by I.N.D.A.M-G.N.A.M.P.A. 2019 and the ``RUDN University Program 5-100''.

\end{document}